\documentclass[11pt]{article}

\usepackage[leqno]{amsmath}
\usepackage{amssymb,amsthm,tabls,url}

\newtheorem{theorem}{Theorem}[section]
\newtheorem{cor}[theorem]{Corollary}
\newtheorem{prop}[theorem]{Proposition}

\newtheoremstyle{efronremark}% entire heading is optional argument []
{6pt}{6pt}{}{}{\itshape}{\quad}{ }{\thmnote{#3}}

\theoremstyle{efronremark}\newtheorem*{eremark}{}

\usepackage[pdftex]{graphicx}

\usepackage{varioref}
\labelformat{section}{Section~#1}
\labelformat{subsection}{Section~#1}
\labelformat{table}{Table~#1}

\usepackage[numbers,sort]{natbib}
\date{February 22, 2011}\numberwithin{equation}{section}

\newcommand\pp{\mathbb{P}}

\begin{document}

\title{Foulkes Characters, Eulerian Idempotents,\\ and an Amazing Matrix}

\author{Persi Diaconis\footnote{Supported in part by NSF grant 0804324.}\\
\textit{Department of Mathematics}\\\textit{Stanford University}\and%
        Jason Fulman\footnote{Supported in part by NSF grant 0802082 and NSA grant H98230-08-1-0133.}\ %
\footnote{Corresponding author: \texttt{fulman@usc.edu}}\\
\textit{Department of Mathematics}\\\textit{University of Southern California}}

\maketitle

\begin{abstract}
  John Holte \cite{holte} introduced a family of ``amazing matrices''
  which give the transition probabilities of ``carries'' when adding a
  list of numbers. It was subsequently shown that these same matrices
  arise in the combinatorics of the Veronese embedding of commutative
  algebra \cite{brenti,pd173,pd186} and in the analysis of riffle
  shuffling \cite{pd173,pd186}. We find that the left eigenvectors of
  these matrices form the \textit{Foulkes character table} of the
  symmetric group and the right eigenvectors are the \textit{Eulerian
    idempotents} introduced by Loday \cite{loday} in work on
  Hochschild homology.  The connections give new closed formulae for
  Foulkes characters and allow explicit computation of natural
  correlation functions in the original carries problem.

  \begin{eremark}[Keywords:]
  Foulkes character, carry, Eulerian idempotent, symmetric group
  \end{eremark}

  \begin{eremark}[AMS 2010 subject classifications:]
    primary 20C30;  secondary 60C05, 60J10.
  \end{eremark}
\end{abstract}

\section{Introduction}\label{sec1}

When $n$ numbers are added in the usual way, ``carries'' accrue along
the way. For example, working base $b=10$, the display shows the
carries along the top when $n=3$ ten-digit numbers are added:
\begin{center}\begin{tabular}{rr}
\it{2} \it{22111} \it{12120}\\
78667\ 51918\\
65921\ 47787\\
88424\ 99859\\\hline
2\ 33013\ 99564
\end{tabular}\end{center}

\noindent Here the carries (reading right to left in bold print) are $\kappa_0=0,\
\kappa_1=2,\ \kappa_2=1,\ \kappa_3=2,\dots$. When $n$ numbers are
added, the carries can be $0,1,2,\dots,n-1$. If the digits are chosen
uniformly at random in $\{0,1,\dots,b-1\}$, it is easy to see that the
carries form a Markov chain: the chance that the next carry is $j$
given the past carries only depends on the last carry. Thus the
distribution of carries is determined by the transition matrix
\begin{equation*}
M(i,j)=\text{chance \{next carry is $j\ |$ last carry is $i$\}}.
\end{equation*}
The carries process was studied by Holte \cite{holte} who showed
\begin{equation}
M(i,j)=\frac1{b^n}\sum_{l=0}^{j-\lfloor i/b\rfloor}(-1)^l\binom{n+1}{l}\binom{n-1-i+(j+1-l)b}{n},
\qquad 0\leq i,j\leq n-1.
\label{11}
\end{equation}
For example, when $n=3$, the matrix is
\begin{equation*}
\frac1{6b^2}\begin{bmatrix}b^2+3b+2&4b^2-4&b^2-3b+2\\b^2-1&4b^2+2&b^2-1\\b^2-3b+2&4b^2-4&b^2+3b+2\end{bmatrix}.
\end{equation*}
Holte found the eigenvalues, eigenvectors, and many amazing properties
of these matrices.

Work of \cite{pd173,pd186,brenti} shows that the same matrix arises in
the analysis of the Gilbert--Shannon--Reeds method of shuffling cards
and in the Hilbert series of the Veronese embedding of projective
varieties.

The main results of the present paper identify a different area where
the matrix appears. The left eigenvectors of the matrix are the
\textit{Foulkes characters} of the symmetric group. The right
eigenvectors are the \textit{Eulerian idempotents} that occur in the
study of free Lie algebras and Hochschild homology. We obtain new
closed-form expressions for these characters.

\ref{sec2} gives background on Foulkes characters and presents some new
results for left eigenvectors. \ref{sec3} does the same for the right
eigenvectors and applies some of the new formulae to the original
carries process, giving the variance and covariance of the number of carries.
\ref{sec4} gives another connection between representation theory of
the symmetric group (the RSK correspondence) and carries.

\section{Foulkes characters}\label{sec2}

This section introduces the Foulkes characters of the symmetric group
and some of their properties (\ref{sec21}). It shows that the Foulkes
characters are the left eigenvectors of the transition matrix $M$ of
\eqref{11} (\ref{sec22}). This connection is used to prove a branching
rule (from $S_n$ to $S_{n-1}$) and a closed-form formula for Foulkes
characters (\ref{sec23}).

\subsection{Background on Foulkes characters}\label{sec21}

Foulkes characters were discovered by \citeauthor{foulkes}
\cite{foulkes} as part of the study of the descent patterns in the
permutation group. They are developed in \cite{kerber84} and
\cite{kerber99} gives a readable textbook treatment.
\citeauthor{gessel} \cite{gessel} use them to enumerate permutations
by descents and conjugacy classes; see \cite{pd77} for a probabilistic
interpretation of these results. \citeauthor{stan07} \cite{stan07}
uses Foulkes characters to develop enumerative results for alternating
permutations by cycle type.

Recall that a permutation $\sigma\in S_n$ has a descent at $i$ if
$\sigma(i+1)<\sigma(i)$. The set of places where descents occur is
$D(\sigma)\subseteq[n-1]$. For example, if
$\sigma=4\underbar{5}123\underbar{7}6,\ D(\sigma)=\{2,6\}$. If
$U\subseteq[n-1]$ is any set, Foulkes suggested constructing a ribbon
shape (also called a rim hook) $R(U)$ beginning with a single box and
sequentially adding the next box below the last box if $i\in U$, and
to the left of the last box if $i\notin U,\ 1\leq i\leq n-1$. Thus, if
$U=\{2,6\}$, boxes are built up as follows:
\begin{figure}[h]
\centering
\includegraphics[scale=0.85, trim=0in 0in 0in 0in, clip=true]{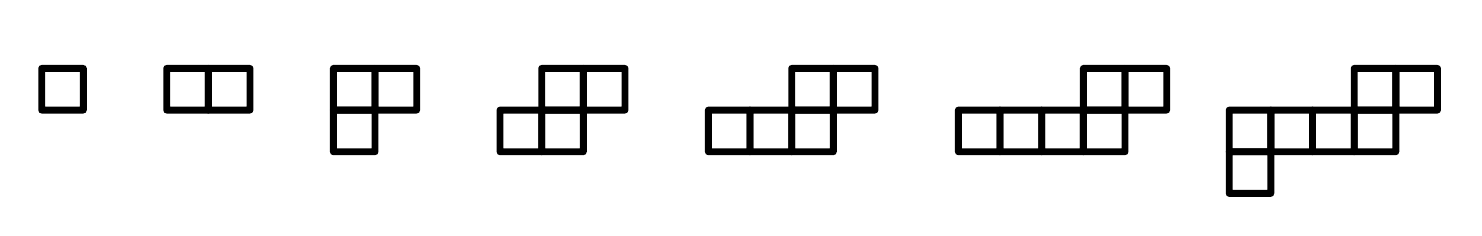}
\end{figure}

\noindent The final skew shape will have $n$ boxes and be the lower
rim of a partition $\alpha$; in this example, $\alpha$ is 5,4,1, and
the ribbon shape is 5,4,1$\backslash$3:
\begin{figure}[h]
\centering
\includegraphics[scale=1.0, trim=0in 0in 0in 0in, clip=true]{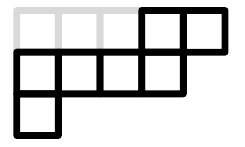}
\end{figure}

\noindent Labeling the boxes in the ribbon shape by all ways they can
be sequentially removed from $\alpha$ and reading this from right to
left and top to bottom gives all permutations with the original $U$ as
descent set. For example, removing boxes in the order shown as
\begin{figure}[h]
\centering
\includegraphics[scale=1.0, trim=0in 0in 0in 0in, clip=true]{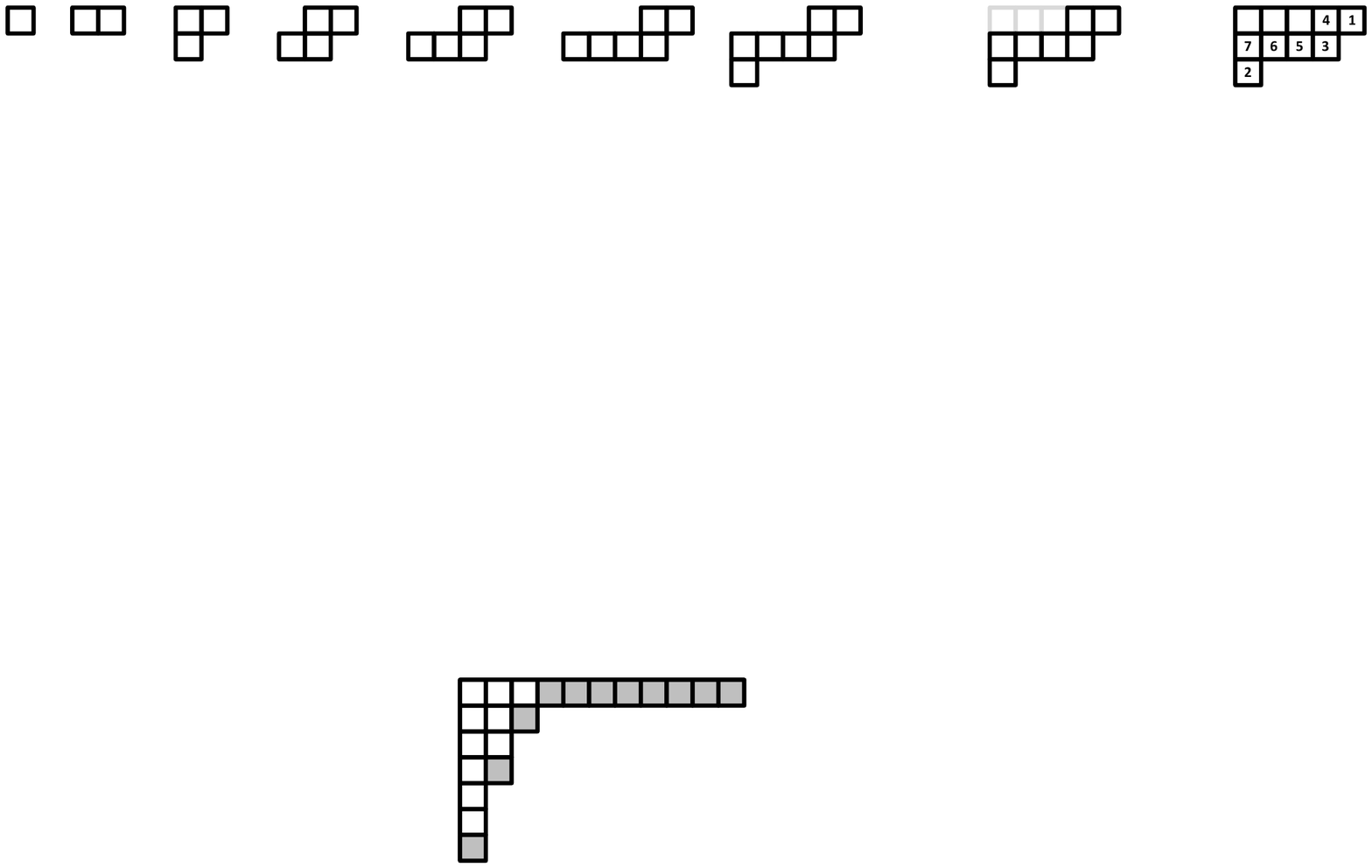}
\end{figure}

\noindent gives 1\underbar{4}356\underbar{7}2. The skew shape $R(U)$
corresponding to $U\subseteq[n-1]$ gives a skew character
$\chi^{R(U)}$: if $R(U)=\alpha\backslash\beta$ and $\chi^\lambda$ is
an irreducible character of the symmetric group $S_n$, the coefficient
of $\chi^\lambda$ in $\chi^{R(U)}$ is
$\langle\chi^\beta\cdot\chi^\lambda|\chi^\alpha\rangle$ (see
\cite[Sect.\ 1.7]{macd}). From the development above, the dimension of
$\chi^{R(U)}$ is the number of permutations with descent set $U$.
Solomon \cite[Sect.\ 6]{solomon} describes a related construction of
MacMahon in his work on Simon Newcomb's problem.
\begin{table}[htb]
\begin{center}\begin{tabular}{r|rrrrr}
&5&4&3&2&1\\\hline
0&1&$-1$&1&$-1$&1\\
1&26&$-10$&2&2&$-4$\\
2&66&0&$-6$&0&6\\
3&26&10&2&$-2$&$-4$\\
4&1&1&1&1&1
\end{tabular}\end{center}
\caption{Foulkes character table for $n=5$.}
\label{table1}
\end{table}

For fixed $k,\ 0\leq k\leq n-1$, the Foulkes character $\chi^{n,k}$ is
defined as the sum of $\chi^{R(U)}$ over all $U$ with $n-k-1$ descents. It
follows that the dimension of $\chi^{n,k}$ is the \textit{Eulerian
  number} $A(n,k)$, the number of permutations with $k$ descents.
Foulkes showed that $\chi^{n,k}(\sigma)$ only depends on $\sigma$
through the number of cycles in $\sigma$. In particular,
\begin{align}
&\chi^{n,n-1}\text{ is the trivial character;}\label{21}\\
&\chi^{n,0}\text{ is the alternating character.}\label{22}
\end{align}
Most importantly, letting $\chi_j^{n,k}$ denote the value of the
Foulkes character on permutations with $j$ cycles (so the dimension
$\chi_n^{n,k}=A(n,k)$),
\begin{equation}
\chi_j^{n,k}=\chi_j^{n-1,k-1}-\chi_j^{n-1,k}\qquad\text{if }j<n,\ k>0.
\label{23}
\end{equation}
This, with the starting value $\chi_1^{1,0}=1$, gives an efficient way
to build a Foulkes character table. Let $k=0,1,\dots,n-1$ index the
rows and $j=n,n-1,\dots,1$ index the columns. \ref{table1} gives the
example when $n=5$.

Further properties of Foulkes characters appear in
\citeauthor{kerber84} \cite{kerber84}:
\begin{align}
&\langle\chi^{n,k},\chi^\lambda\rangle>0\Longrightarrow\lambda_1\leq
k+1,\ \lambda_1'\leq n-k.\label{24}\\
&\langle\chi^{n,k},\chi^{j+1,1^{n-j-1}}\rangle>0\Leftrightarrow j=k.\label{25}\\
&\text{The $\chi^{n,k}$ are linearly independent.}\label{26}\\
&\text{If $\chi:S_n\to\mathbb{R}$ denotes a character, depending only on
  the number of cycles, then}\label{27}\\
&\qquad\chi=\sum_i\frac{\langle\chi,\chi^{i+1,1^{n-i-1}}\rangle}{\chi^{i+1,1^{n-i-1}}(id)}\chi^{n,i}.\notag
\end{align}
Thus the $\chi^{n,k}$ form a $\mathbb{Q}$ basis for the characters
that only depends on the number of cycles. Hidden in the proof of
\eqref{27}: the hook character $\chi^{i+1,1^{n-i-1}}$ is the only hook
occuring in $\chi^{n,i}$ and it occurs with multiplicity its degree
$\binom{n-1}{i}$. A related fact appears in \citeauthor{solomon}
\cite[Th.\ 4]{solomon}. \citeauthor{kerber84} \cite{kerber84} further
determine the permutation character for $S_n$ acting on $[M]^n$:
\begin{equation}
\chi_M(\sigma):=M^{\text{\# cycles in }(\sigma)}
\label{28}
\end{equation}
has the decomposition
\begin{equation*}
\chi_M=\sum_k\binom{M+k}{n}\chi^{n,k}.
\end{equation*}

The $\chi^{n,k}$ are usually \textit{not} irreducible, and
\citeauthor{kerber84} \cite{kerber84} give an interesting
combinatorial rule for decomposing $\chi^{R(U)}$ (and thus
$\chi^{n,k}$). They show that $\chi^{n,k}$ is a sum of Lefschetz
characters for $S_n$ acting on $2^{[n]}$.

Marty Isaacs conjectured that $n!$ divides the determinant of
the Foulkes character table. In fact the following is true:
\begin{equation}
\text{The determinant of the $n\times n$ matrix with $k,j$ entry
  $\chi_j^{n,k}$ is $n!(n-1)! \cdots 2!$.}
\label{new29}
\end{equation}
\begin{proof}
  Construct an $n\times n$ matrix $A$ from the $(n-1)\times(n-1)$
  Foulkes character table by adding a left column consisting of the
  partial sums of the Eulerian numbers $A(n,0),A(n,0)+A(n,1),\dots,n!$
  and filling out the rest of the top row with zeros. Thus, when
  $n=5$,
\begin{equation*}
A=\begin{pmatrix}
120&0&0&0&0\\119&1&-1&1&-1\\93&11&-3&-1&3\\27&11&3&-1&-3\\1&1&1&1&1\end{pmatrix}.
\end{equation*}
The $4\times4$ matrix in the lower right corner is the Foulkes
character table for $n=4$. The first column entries are the partial
sums of the Eulerian numbers 1, 26, 66, 26, 1. In particular, the
(1,1) entry is $n!$, so by induction the determinant of $A$ is
$n!(n-1)! \cdots 2!$.

The $n\times n$ Foulkes character table is constructed from $A$ as
follows: in $A$, subtract row 2 from row 1, then row 3 from row 2, and
so on. The recurrence \eqref{23} shows this gives the $n\times n$
Foulkes character table.
\end{proof}

Marty Isaacs observes that while the $\{\chi^{n,k}\}$ are not
disjoint, they sum up to the regular character of $S_n$; see
\cite[Th.\ 2]{solomon} for a proof.  Alas, this does not
seem to be enough to have the nice theory of supercharacters
\cite{aguiar} carry over, but the parallels are intriguing. Further
properties of Foulkes characters are given in \ref{sec23} after the
connection with the carries transition matrix is developed.

Rim hook characters are a basic construction of representation theory
of $S_n$; see \cite{billera,lascoux} and the references there. They are also
available for other Coxeter groups \cite{solomon}. Foulkes' innovation,
showing that sums of these characters have interesting properties, has
not been explored for general type.

\subsection{The connection with carries}\label{sec22}

We noticed from Holte's paper \cite{holte} that the matrix of left
(row) eigenvectors for the carries Markov chain on $\{0,1,2,3,4\}$
(e.g., working base 5) is
\begin{equation*}\begin{pmatrix}
1&26&66&26&1\\1&10&0&-10&-1\\1&2&-6&2&1\\1&-2&0&2&-1\\1&-4&6&-4&1\end{pmatrix}.
\end{equation*}
Comparing this with \ref{table1}, the Foulkes character table, leads
to the following result.
\begin{theorem}\label{leigen}
  Let $v^n_{i,j}$ denote the $j$th entry of the $i$th left eigenvector
  of the carries matrix for addition of $n$ numbers base $b$ (here
  $0\leq i,j\leq n-1$, and the eigenvalues are $1/b^i$). Then
\begin{equation*}
v^n_{i,j}=\chi^{n,n-j-1}_{n-i}.
\end{equation*}
\end{theorem}
\begin{proof}
  The first case is that $i=0$. From \cite{holte}, $v^n_{0,j}=A(n,j)$.
  From the dimension formula, $\chi^{n,n-j-1}_n=A(n,n-1-j)$. By
  symmetry of the Eulerian numbers, $A(n,j)=A(n,n-1-j)$, so the
  theorem follows in the first case.

  The second case is that $j=n-1$. By equation \eqref{22},
  $\chi^{n,0}_j=(-1)^{n-j}$. Thus we need to show that
  $v^n_{n-j,n-1}=(-1)^{n-j}$. By Holte's formula for the left
  eigenvectors of the carries chain \cite[p.\ 143]{holte}, it follows
  that
\begin{equation*}
v^n_{n-j,n-1}=\sum_{r=0}^{n-1} (-1)^r \binom{n+1}{r} (n-r)^j.
\end{equation*}
The result now follows by induction, since
\begin{align*}
&\sum_{r=0}^{n-1}(-1)^r\binom{n+1}{r}(n-r)^{j-1}(n-r)\\
&\qquad=n\sum_{r=0}^{n-1}(-1)^r\binom{n+1}{r}(n-r)^{j-1}-\sum_{r=0}^{n-1}(-1)^rr\binom{n+1}{r}(n-r)^{j-1}\\
&\qquad=n(-1)^{n-j+1}-(n+1)\sum_{r=1}^{n-1}(-1)^r\binom{n}{r-1}(n-r)^{j-1}\\
&\qquad=n(-1)^{n-j+1}+(n+1)\sum_{r=0}^{n-2}(-1)^r\binom{n}{r}(n-1-r)^{j-1}\\
&\qquad=n(-1)^{n-j+1}+(n+1)(-1)^{n-j}\\
&\qquad=(-1)^{n-j}.
\end{align*}

For the remaining cases, $i>0$ and $j<n-1$. By the recursive formula
\eqref{23}, it is enough to show that
\begin{equation*}
v^n_{i,j} = v^{n-1}_{i-1,j} - v^{n-1}_{i-1,j-1}
\end{equation*}
for $i>0,\ j<n-1$. From page 144 of \cite{holte},
\begin{align*}
v^n_{i,j}=&\text{ Coefficient of $x^{j+1}$ in $(1-x)^{n+1} \left(x\frac{d}{dx}\right)^{n-i}(1-x)^{-1}$}.\\
%\end{equation*}
\intertext{Clearly}
%\begin{align*}
&\text{Coefficient of $x^{j+1}$ in $(1-x)^{n+1} \left(x \frac{d}{dx} \right)^{n-i}(1-x)^{-1}$}\\
&\qquad=\text{ Coefficient of $x^{j+1}$ in $(1-x)^{n} \left(x \frac{d}{dx} \right)^{n-i}(1-x)^{-1}$}\\
 &\qquad\qquad -\text{ Coefficient of $x^{j}$ in $(1-x)^{n} \left(x \frac{d}{dx} \right)^{n-i}(1-x)^{-1}$}
\end{align*}
which implies the result.
\end{proof}

\subsection{Some consequences}\label{sec23}

In \cite{holte}, Holte gave a closed formula for the left
eigenfunctions,
\begin{equation}
v_{i,j}^n=\sum_{r=0}^{j+1}(-1)^r\binom{n+1}{r}(j+1-r)^{n-i}.
\label{210}
\end{equation}
Thus we get an apparently new formula for the Foulkes characters.
\begin{cor}
\begin{equation*}
\chi_j^{n,k}=\sum_{r=0}^{n-k}(-1)^r\binom{n+1}{r}(n-k-r)^j.
\end{equation*}
\label{cor12}
\end{cor}
In rereading Foulkes \cite[Sect.\ 4]{foulkes} we found the formula
(Th.\ 4.1),
%of Theorem \ref{desgen},
\begin{equation*}
\chi_j^{n,k}=\sum_{r=0}^{n-j}(-1)^r\binom{n-j}{r}A(j,k+j+r-n).
\end{equation*}
This seems a little less direct than Corollary \ref{cor12}. Corollary
\ref{cor12} gives a direct proof of the following restriction formula
of Foulkes characters from $S_n$ to $S_{n-1}$.
\begin{cor}[{\cite[Cor.\ 4.6]{foulkes}}]
\begin{equation*}
\chi_{S_{n-1}}^{n,k}=(k+1)\chi^{n-1,k}+(n-k)\chi^{n-1,k-1},\qquad0\leq k\leq n-1.
\end{equation*}
\label{cor13}
\end{cor}
\begin{eremark}[Remarks]
  We first learned Corollary \ref{cor13} from Marty Isaacs, who both
  observed it and showed that it follows from Corollary \ref{cor12}.
  The well-known recursion formula for the Eulerian numbers
  $A(n,k)=(k+1)A(n-1,k)+(n-k)A(n-1,k-1)$ is the special case of
  evaluation at the identity. Thus Corollary \ref{cor13} represents a
  ``categorification'' of this recurrence.
\end{eremark}
\begin{proof}
The required formula translates to
\begin{equation*}
\chi_j^{n,k}=(k+1)\chi_{j-1}^{n-1,k}+(n-k)\chi_{j-1}^{n-1,k-1},\qquad0\leq k\leq n-1,\ j\geq1.
\end{equation*}
From Corollary \ref{cor12},
\begin{align*}
(k+1) \chi^{n-1,k}_{j-1} &+ (n-k) \chi^{n-1,k-1}_{j-1}\\
& =  (k+1) \sum_{r=0}^{n-k-1} (-1)^r \binom{n}{r} (n-k-r-1)^{j-1} + (n-k) \chi^{n-1,k-1}_{j-1} \\
& =  (k+1) \sum_{r=1}^{n-k} (-1)^{r-1} \binom{n}{r-1} (n-k-r)^{j-1} \\
&\qquad + (n-k)^j  + (n-k) \sum_{r=1}^{n-k} (-1)^r \binom{n}{r} (n-k-r)^{j-1} \\
& = (n-k)^j + \sum_{r=1}^{n-k} (-1)^{r-1} \binom{n}{r-1} (n-k-r)^{j-1} \\
&\qquad + k \sum_{r=1}^{n-k} (-1)^{r-1} \binom{n}{r-1} (n-k-r)^{j-1} \\
&\qquad + n \sum_{r=1}^{n-k} (-1)^r \binom{n}{r} (n-k-r)^{j-1} - k \sum_{r=1}^{n-k} (-1)^r \binom{n}{r} (n-k-r)^{j-1}\\
\end{align*}
\begin{align*}
& = (n-k)^j + n \sum_{r=1}^{n-k} (-1)^r \binom{n}{r} (n-k-r)^{j-1} \\
&\qquad + \sum_{r=1}^{n-k} (-1)^{r-1} \binom{n}{r-1} (n-k-r)^{j-1} \\
&\qquad - k \sum_{r=1}^{n-k} (-1)^r \binom{n+1}{r} (n-k-r)^{j-1} - \sum_{r=1}^{n-k} (-1)^r r \binom{n+1}{r} (n-k-r)^{j-1} \\
&\qquad + \sum_{r=1}^{n-k} (-1)^r r \binom{n+1}{r} (n-k-r)^{j-1} + n \sum_{r=1}^{n-k} (-1)^r \binom{n+1}{r} (n-k-r)^{j-1} \\
&\qquad - n \sum_{r=1}^{n-k} (-1)^r \binom{n+1}{r} (n-k-r)^{j-1} \\
& = \left[ (n-k)^j + \sum_{r=1}^{n-k} (-1)^r \binom{n+1}{r} (n-k-r)^j \right] \\
&\qquad + n \sum_{r=1}^{n-k} (-1)^r \binom{n}{r} (n-k-r)^{j-1} + \sum_{r=1}^{n-k} (-1)^{r-1} \binom{n}{r-1} (n-k-r)^{j-1} \\
&\qquad - \sum_{r=1}^{n-k} (-1)^r (n-r) \binom{n+1}{r} (n-k-r)^{j-1} \\
& = \chi^{n,k}_j + \sum_{r=1}^{n-k} (-1)^r (n-k-r)^{j-1} \left[ n \binom{n}{r} - \binom{n}{r-1} - (n-r) \binom{n+1}{r} \right] \\
& = \chi^{n,k}_j
\end{align*}
where the final equality used the identity
\begin{equation*}
n\binom{n}{r}-\binom{n}{r-1}-(n-r)\binom{n+1}{r}=0.\qedhere
\end{equation*}
\end{proof}

\section{Riffle shuffle idempotents and right eigenvectors of the carries matrix}\label{sec3}

This section simplifies Holte's formula for the right eigenvectors of
the carries matrix, and relates these eigenvectors to representation
theory of the symmetric group (\ref{sec31}). The eigenfunctions are
used to compute basic things about carries in \ref{sec32}.

\subsection{Right Eigenfunctions}\label{sec31}

To begin, let $u^n_j(i)$ denote the value of the $j$th right
eigenvector of the carries chain ($0 \leq j \leq n-1$, eigenvalues
$1/b^j$) evaluated at $i$ ($0 \leq i \leq n-1$). We also let $s(n,k)$
be the Stirling number of the first kind, defined as $(-1)^{n-k}$
multiplied by the number of permutations on $n$ symbols with $k$
cycles. It can also be defined by the equation
\begin{equation}\label{stirdef}
x(x-1)\dots(x-n+1)=\sum_{k\geq0}s(n,k)x^k.
\end{equation}
Theorem 4 of Holte \cite{holte} shows that
\begin{equation*}
u^n_j(i) = \sum_{k=n-j}^n s(n,k) \binom{k}{n-j} (n-1-i)^{k-(n-j)}
\end{equation*}
where $0^0$ is taken to be 1. Note that $u_j^n(i)$ is a polynomial in
$i$ of degree $j$. For example, $u_0^n(i)=1,\
u_1^n(i)=n(n-1-i)-\binom{n}2$.

The next theorem gives a simpler formula for $u^n_j(i)$.
\begin{theorem}\label{simpler}
\begin{align*}
u^n_j(i)&=n!\sum_{k\geq0}\frac{s(k,n-j)}{k!}\binom{n-i-1}{n-k}\\
&=n!\cdot\textup{Coefficient of $x^{n-j}$ in $\binom{x+n-i-1}{n}$}.
\end{align*}
\end{theorem}
\begin{proof}
  By Theorem 4 of Holte \cite{holte}, $u^n_j(i)$ is $n!$ multiplied by
  the the $i,j$ entry of the inverse of the matrix of left row
  eigenvectors. From this and equation \eqref{210}, proving the
  first equality of the theorem is equivalent to proving that
\begin{equation*}
n!\sum_{l\geq0}\sum_{k\geq0}\frac{s(k,n-l)}{k!}\binom{n-i-1}{n-k}\sum_{r=0}^{j+1}(-1)^r\binom{n+1}{r}(j+1-r)^{n-l}
 =n!\delta_{i,j}.
\end{equation*}
Now
\begin{align*}
&n!\sum_{l=n-k}^n\sum_{k \geq 0}\frac{s(k,n-l)}{k!}\binom{n-i-1}{n-k}\sum_{r=0}^{j+1}(-1)^r\binom{n+1}{r}(j+1-r)^{n-l}\\
&\qquad =n!\sum_{k\geq0}\binom{n-i-1}{n-k}\sum_{r=0}^{j+1}(-1)^r\binom{n+1}{r}\frac1{k!}\sum_{l=n-k}^ns(k,n-l)(j+1-r)^{n-l}\\
&\qquad =n!\sum_{k\geq0}\binom{n-i-1}{n-k}\sum_{r=0}^{j+1}(-1)^r\binom{n+1}{r}\frac1{k!}\sum_{t=0}^ks(k,t)(j+1-r)^{t}\\
&\qquad =n!\sum_{k\geq0}\binom{n-i-1}{n-k}\sum_{r\geq0}(-1)^r\binom{n+1}{r}\binom{j+1-r}{k}\\
&\qquad =n!\sum_{r\geq0}(-1)^r\binom{n+1}{r}\sum_{k\geq0}\binom{n-i-1}{n-k}\binom{j+1-r}{k}\\
&\qquad =n!\sum_{r\geq0}(-1)^r\binom{n+1}{r}\binom{n+j-i-r}{n}\\
&\qquad =n!\delta_{i,j}.
\end{align*}
The third equality used equation \eqref{stirdef}, the fifth equality
used the basic identity
$\sum_i\binom{a}{i}\binom{b}{n-i}=\binom{a+b}{n}$ \cite[p.\
12]{stan97}, and the final equality is from page 147 of \cite{holte}.

To prove the second equality of the theorem, write
\begin{equation*}
\binom{x+n-i-1}{n}=\sum_k\binom{x}{k}\binom{n-i-1}{n-k}
\end{equation*}
and use equation \eqref{stirdef}.
\end{proof}
\begin{eremark}[Remark]
Let $E_{n,k}$ be elements of the symmetric group algebra defined by
the equation
\begin{equation*}
\sum_{k=1}^nx^kE_{n,k}=\sum_{w\in S_n}\binom{x+n-d(w)-1}{n}
\end{equation*}
where $d(w)$ denotes the number of descents of $w$. By work of
\citeauthor{garsia89} \cite{garsia89}, these are orthogonal
idempotents of the symmetric group algebra whose sum is the identity.
They also arise in the theory of riffle shuffling \cite{pd92} and in
Hochschild homology \cite{hanlon90}. Their images under the sign map
are known as \textit{Eulerian idempotents}. In this version, they were
discovered by Gerstenhaber--Schack \cite{gersten} to give Hodge
decompositions of Hochschild homology. They have been developed by
Loday \cite{loday} for cyclic homology. Patras \cite{patras} gives an
unusual treatment involving decompositions of the $n$-cube into
simplices. For a textbook treatment, see Weibel \cite[Sect.\
9.4.3]{weibel}.

The eigenvectors of the carries and descent matrix lift to eigenvectors
of the full riffle shuffle matrix. These in turn are identified in
Denham \cite{denham} and Diaconis--Ram \cite{pd190}.
\end{eremark}

Clearly the value of $E_{n,k}$ on a permutation depends only on its
number of descents.  Letting $E_{n,k}(d)$ denote the value of
$E_{n,k}$ on a permutation with $d$ descents, we have the following
corollary of Theorem \ref{simpler}.
\begin{cor}\label{val}
\begin{equation*}
u^n_j(i) = n! E_{n,n-j}(i).
\end{equation*}
\end{cor}

It would be nice to have a more conceptual proof of Corollary
\ref{val}.

\subsection{Applications}\label{sec32}

This section gives some applications of the explicit form of the right
eigenvectors of the carries chain for the addition of $n$ numbers base
$b$. We note that another application (to lower bounding the
convergence rate of the carries chain) appears in \cite{pd186}.

The transition matrix of the carries chain is viewed as a linear
operator on functions in the usual way: $Kf[x]=\sum_y K(x,y) f(y)$.
Let $\kappa_r$ denote the value of the carry from column $r-1$ to
column $r$.
\begin{prop}\label{stat}
  Suppose that the carries chain is started from its stationary
  distribution $\pi$. Then for $n \geq 2$, the covariance
\begin{equation*}
\textup{Cov}(\kappa_0\kappa_r)=\frac1{b^r}\frac{n+1}{12}>0.
\end{equation*}
\end{prop}
\begin{proof}
Clearly
\begin{equation*}
E(\kappa_0\kappa_r)=\sum_ii\pi(i)E(\kappa_r|\kappa_0=i).
\end{equation*}
By Theorem \ref{simpler}, $u^n_1$ is a right eigenvector of the
carries chain with eigenvalue $1/b$. Dividing by $-1/n$, we have the
eigenvector $f(i)=i-\frac{n-1}2$,
\begin{equation*}
E(\kappa_r|\kappa_0=i)=K^r\left(f+\frac{n-1}2\right)[i]=\frac1{b^r}\left(i-\frac{n-1}2\right)+\frac{n-1}2.
\end{equation*}
It follows that
\begin{equation}\label{one}
E(\kappa_0\kappa_r)=\sum_ii\left[\frac1{b^r}\left(i-\frac{n-1}2\right)+\frac{n-1}2\right]\pi(i).
\end{equation}

From \cite {holte} or \cite{pd186}, the stationary distribution of the
carries chain is $\pi(i)=A(n,i)/n!$ (here $A(n,i)$ is the number of
permutations on $n$ symbols with $i$ descents), and it is well known
that for $n \geq 2$ the mean and variance of the Eulerian numbers are
$\frac{n-1}{2}$ and $\frac{n+1}{12}$. This, together with equation
\eqref{one}, gives
\begin{align*}
E(\kappa_0\kappa_r)&=
\frac1{b^r}\left[\frac{n+1}{12}+\left(\frac{n-1}2\right)^2\right]
    -\frac1{b^r}\left(\frac{n-1}2\right)^2+\left(\frac{n-1}2\right)^2\\
&=\frac1{b^r}\frac{n+1}{12}+\left(\frac{n-1}2\right)^2
\end{align*}
and the result follows since
$E(\kappa_0)E(\kappa_r)=\left(\frac{n-1}2\right)^2$.
\end{proof}

A similar calculation allows us to compute the covariance started from
the state 0:
\begin{prop}\label{use1}
  Suppose that the carries chain is started from $0$, and that
  $n\geq2$. Then
\begin{enumerate}
\item $E(\kappa_r)=\left(1-\frac1{b^r}\right)\frac{n-1}2$.
\item $\textup{Var}(\kappa_r)=\left(1-\frac1{b^{2r}}\right)\frac{n+1}{12}$.
\item $\textup{Cov}(\kappa_s,\kappa_{s+r})=\frac1{b^r}\frac{n+1}{12}\left(1-\frac1{b^{2s}}\right)>0$.
\end{enumerate}
\end{prop}
\begin{proof}
  Parts 1 and 2 are proved in Theorem 4.1 of \cite{pd173}. For part 3,
  arguing as in the proof of Proposition \ref{stat},
\begin{align*}
E(\kappa_s\kappa_{s+r})&=\sum_iiP(\kappa_s=i)E(\kappa_{s+r}|\kappa_s=i)\\
&=\sum_iiP(\kappa_s=i)\left[\frac1{b^r}\left(i-\frac{n-1}2\right)+\frac{n-1}2\right]\\
&=\frac1{b^r}\left[\text{Var}(\kappa_s)+E(\kappa_s)^2\right]+\left(1-\frac1{b^r}\right)\frac{n-1}2 E[\kappa_s].
\end{align*}
The result follows from parts 1 and 2 after elementary simplifications.
\end{proof}

The eigenvectors of the carries matrix can be used to give a simple
proof that the sum of $\chi^{n,k}$ gives the regular character. Holte
\cite{holte} shows that if $V$ is the matrix whose rows are left
eigenvectors and $U$ is the matrix whose columns are right
eigenvectors, then $VU=n!\times I$. Just looking at products involving
the first column of $U$ (which is identically one) we get that the
$\chi^{n,k}$ sum up to the regular character.

It is instructive to see the problems encountered in trying to use the
available eigenstructure to bound the rate of convergence of the
carries chain to its stationary distribution $A(n,j)/n!$. Let
$M_b(i,j)$ be the transition matrix \eqref{11} corresponding to adding
$n$ numbers base $b$. Let $M_b^k(i,j)$ be the $k$th power of this
matrix. From elementary linear algebra,
\begin{equation*}
M_b^k(i,j)=\frac{A(n,j)}{n!}+\sum_{a=1}^{n-1}\frac1{b^{ak}}r_a(i)l_a(j)
\end{equation*}
where $l_a,r_a$ are the left and right eigenvectors of $M_b$ normed
  so that $l_0(j)=\frac{A(n,j)}{n!}$ (and $r_0(j)=1$) for $0\leq j\leq
  n-1$. Here
\begin{align*}
l_a(j)&=\frac1{n!}\chi_{n-a}^{n,n-j-1}\text{ with
  $\chi_{n-a}^{n,n-j-1}$ the Foulkes characters of \ref{sec2}.}\\
r_a(j)&=u^n_a(j)\text{ with $u^n_a(j)$ given in Theorem \ref{simpler}.}
\end{align*}
The total variation distance to stationarity, starting at $i=0$, after
$k$ steps is
\begin{equation*}
\frac12\sum_{j=0}^{n-1}\left|M_b^k(0,j)-\frac{A(n,j)}{n!}\right|.
\end{equation*}
From the formulae above,
\begin{equation*}
\left|M_b^k(0,j)-\frac{A(n,j)}{n!}\right|=
\left|\sum_{a=1}^{n-1}\frac1{b^{ak}}r_a(0)l_a(j)\right|.
\end{equation*}
While bounding this is feasible, it is a bit of a mess. In
\cite[Sect.\ 3]{pd186}, a different representation is used to prove
that the carries chain is close to stationarity after order
$\frac12\log_b(n)$ steps.

\section{Carries and the RSK correspondence}\label{sec4}

In this section we use the RSK correspondence to derive a generating
function for descents after a $b^r$-riffle shuffle on a deck of $n$
cards. See \cite{pd92} for background on riffle shuffles. By a main
result of \cite{pd173}, this is equal to the generating function for
the rth carry $\kappa_r$ when $n$ numbers are added base $b$ (and one
can give another proof of Theorem \ref{desgen} using Holte's formula
for $P(\kappa_r=i)$).
\begin{theorem}\label{desgen}
  Let $w$ be produced by a $b^r$-riffle shuffle on a deck of $n$
  cards, and let $d(w)$ denote the number of descents of $w$. Then
\begin{equation}\label{des}
\sum_{w\in S_n}\pp(w)t^{d(w)+1}=\frac{(1-t)^{n+1}}{b^{rn}}\sum_{k\geq1}t^k\binom{b^rk+n-1}{n}.
\end{equation}
\end{theorem}
\begin{proof}
  Let $w$ be a permutation produced by a $b^r$ riffle shuffle. The RSK
  correspondence associates to $w$ a pair of standard Young tableaux
  $(P(w),Q(w))$ of the same shape. Moreover, there is a notion of
  descent set for standard Young tableaux, and by Lemma 7.23.1 of
  \cite{stan99}, the descent set of $w$ is equal to the descent set of
  $Q(w)$. It is known from \cite{fulman02-2} that if $w$ is produced
  by a $b^r$ shuffle, then the chance that $Q(w)$ is equal to any
  particular standard Young tableau of shape $\lambda$ is $s_{\lambda}
  \left(\frac1{b^r},\dots,\frac1{b^r}\right)$, where there are $b^r$
  variables. Letting $f_{\lambda}(a)$ denote the number of standard
  Young tableaux of shape $\lambda$ with $a$ descents, it follows that
\begin{equation*}
\pp\left(d(w)=a\right)=\sum_{|\lambda|=n}f_{\lambda}(a)s_{\lambda}\left(\frac1{b^r},\dots,\frac1{b^r}\right).
\end{equation*}

By equation 7.96 of \cite{stan99},
\begin{equation*}
\sum_{a\geq0}f_{\lambda}(a)t^{a+1}=(1-t)^{n+1}\sum_{k\geq1}s_{\lambda}(1,\dots,1)t^k
\end{equation*}
where in the kth summand, $s_{\lambda}(1,\dots,1)$ denotes the Schur
function with $k$ variables specialized to 1. Thus
\begin{align*}
&\sum_{a\geq0}\pp\left(d(w)=a\right)\cdot t^{a+1}\\
&\qquad = \sum_{a\geq0}\sum_{|\lambda|=n}f_{\lambda}(a)s_{\lambda}\left(\frac1{b^r},\dots,\frac1{b^r}\right)\cdot t^{a+1}\\
&\qquad = (1-t)^{n+1}\sum_{k\geq1}t^k\sum_{|\lambda|=n}s_{\lambda}\left(\frac1{b^r},\dots,\frac1{b^r}\right)s_{\lambda}(1,\dots,1)\\
&\qquad = (1-t)^{n+1}\sum_{k\geq1}t^k[u^n]\sum_{n\geq0}\sum_{|\lambda|=n}s_{\lambda}\left(\frac1{b^r},
               \dots,\frac1{b^r}\right)s_{\lambda}(1,\dots,1)\cdot u^n
\end{align*}
where $[u^n] g(u)$ denotes the coefficient of $u^n$ is a power series
$g(u)$. Applying the Cauchy identity for Schur functions \cite[p.\
322]{stan99}, this becomes
\begin{equation*}
(1-t)^{n+1}\sum_{k\geq1}t^k[u^n](1-u/b^r)^{-b^rk}=\frac{(1-t)^{n+1}}{b^{rn}}\sum_{k\geq1}t^k\binom{b^rk+n-1}{n}
\end{equation*}
as desired.
\end{proof}

\section*{Acknowledgement}

We thank Danny Goldstein, Bob Guralnick, Marty Isaacs, Eric Rains, and
Arun Ram for their insights, hard work, and suggestions.

\bibliography{PDall,PDnone}

\end{document}